# Time–space harmonic polynomials relative to a Lévy process

JOSEP LLUÍS SOLÉ* and FREDERIC UTZET**

*Departament de Matemàtiques, Facultat de Ciències, Universitat Autónoma de Barcelona, 08193 Bellaterra (Barcelona), Spain. E-mail: \*jllsole@mat.uab.cat; \*\*utzet@mat.uab.cat*

In this work, we give a closed form and a recurrence relation for a family of time–space harmonic polynomials relative to a Lévy process. We also state the relationship with the Kailath–Segall (orthogonal) polynomials associated to the process.

## 1. Introduction

Given a stochastic process $X = \{X_t, t \in \mathbb{R}_+\}$ with finite moments of convenient order, a *time–space harmonic polynomial relative* to $X$ is a polynomial $Q(x,t)$ such that the process $M_t = Q(X_t, t)$ is a martingale with respect to the filtration associated with $X$. Major examples are the Hermite polynomials relative to a Brownian motion, the Charlier polynomials relative to a Poisson process and the Laguerre polynomials relative to a Gamma process; for these and other examples, see Feinsilver [2], Schoutens and Teugels [11], Schoutens [10] and Barrieu and Schoutens [1].

Finding a family of time–space harmonic polynomials relative to a particular process is not an easy task; see, for example, Schoutens [10]. Nevertheless, when we are dealing with a Lévy process $X$ that has moment generating function in a neighborhood of the origin, Sengupta [13] mentions following Neveu [7] a general procedure based on the associated exponential martingale (the left-hand side of (1) below). Specifically, that martingale is an analytic function (of $u$) in some neighborhood of the origin and the corresponding Taylor expansion

$$\frac{\exp\{uX_t\}}{\mathbb{E}[\exp\{uX_t\}]} = \sum_{n=0}^{\infty} R_n(X_t, t) \frac{u^n}{n!} \quad (1)$$

has coefficients such that $R_n(x,t)$ are time–space harmonic polynomials. However, when the process does not have a finite moment generating function, this approach does not work; then, Sengupta [13] uses a discretization procedure to extend the results proven by Goswami and Sengupta [3], without explicitly showing the polynomials.

In this paper, using Itô's formula in our proof, we give a closed form for a family of time–space harmonic polynomials of a Lévy process with finite moments of all orders. From an interesting property observed by Sengupta [13], it is deduced that these polynomials play the role of the basic time–space harmonic polynomials because they span the set of all time–space harmonic polynomials. These results provide new insight into time-space harmonic polynomials and offer a general method for computing such polynomials. Furthermore, the construction has a direct extension to the case where the underlying Lévy process has only a finite number of finite moments and we are only interested in a few space–time harmonic polynomials. We also present some examples and study the relationship between time–space harmonic polynomials and the orthogonal polynomials called Kailath–Segall polynomials (see Segall and Kailath [12] or Meyer [6]).

The time–space harmonic polynomials that we propose here are related to the polynomials that give the moments of a random variable as functions of the cumulants. So we will start with some results from the elementary theory of cumulants.

## 2. Moments and cumulants

Consider a random variable $Y$ with finite moments of all orders. Denote by $\mu_n$ (resp. $\kappa_n$) its moment (resp. its cumulant) of order $n$. It is well known that $\mu_n$ can be written as a polynomial in $\kappa_1, \ldots, \kappa_n$, for example,

$$\mu_1 = \kappa_1, \qquad \mu_2 = \kappa_1^2 + \kappa_2, \qquad \mu_3 = \kappa_1^3 + 3\kappa_1\kappa_2 + \kappa_3, \ldots .$$

If $Y$ has moment generating function in some open interval containing 0, a general expression for $\mu_n$ can be deduced from the relationship between the moment generating function and the cumulant generating function:

$$\exp\left\{\sum_{n=1}^{\infty} \kappa_n \frac{u^n}{n!}\right\} = \sum_{n=0}^{\infty} \mu_n \frac{u^n}{n!}. \qquad (2)$$

When $Y$ does not fulfill the condition in the moment generating function, we can consider the finite series up to the convenient order; see McCullagh [5], page 26, for the justification of this procedure.

Through a direct examination of (2), Kendall and Stuart ([4], formula (3.33)) deduce that

$$\mu_n = n! \sum_{m=1}^{n} \sum \frac{\kappa_1^{r_1}}{(1!)^{r_1} r_1!} \cdots \frac{\kappa_m^{r_m}}{(m!)^{r_m} r_m!}, \qquad (3)$$

where the second summation is over all non-negative integers $r_1, \ldots, r_m$ such that $\sum_{j=1}^{m} j r_j = n$.

Let $\Gamma_n(x_1, \ldots, x_n)$, $n \geq 1$, be the polynomial defined by the expression (3). That is, we have

$$\mu_n = \Gamma_n(\kappa_1, \ldots, \kappa_n). \tag{4}$$

Also, write $\Gamma_0 = 1$.

A recurrence formula follows from (2) for $\Gamma_n$ ([15], Proposition 5.1.7):

$$\Gamma_{n+1}(x_1, \ldots, x_{n+1}) = \sum_{j=0}^{n} \binom{n}{j} \Gamma_j(x_1, \ldots, x_j) \, x_{n+1-j}. \tag{5}$$

Further,

$$\frac{\partial \Gamma_n(x_1, \ldots, x_n)}{\partial x_j} = \binom{n}{j} \Gamma_{n-j}(x_1, \ldots, x_{n-j}), \qquad j = 1, \ldots, n. \tag{6}$$

From the Taylor expansion of $\Gamma_n(x_1 + y, x_2, \ldots, x_n)$ at $y = 0$, we obtain the following expression that we will need later:

$$\Gamma_n(x_1 + y, x_2, \ldots, x_n) = \sum_{j=0}^{n} \binom{n}{j} \Gamma_{n-j}(x_1, \ldots, x_{n-j}) \, y^j. \tag{7}$$

Interchanging the roles of $x_1$ and $y$ and evaluating the function at 0, we obtain

$$\Gamma_n(x_1, x_2, \ldots, x_n) = \sum_{j=0}^{n} \binom{n}{j} \Gamma_{n-j}(0, x_2, \ldots, x_{n-j}) \, x_1^j. \tag{8}$$

Applying this formula to $\kappa_1, \ldots, \kappa_n$ and using the fact that $\kappa_2, \ldots, \kappa_n$ are the cumulants of both $Y$ and $Y - \mathbb{E}Y$ gives the expression for the non-centered moment $\mu_n = \mathbb{E}[Y^n]$ (on the left) from the centered moments $\mathbb{E}[(Y - \mathbb{E}Y)^{n-j}]$ and the powers of the expectation $(\mathbb{E}Y)^j$. It is a classical formula (see [4], formula (3.6)).

## 3. Time–space harmonic polynomials for a Lévy process

Let $X = \{X_t, t \geq 0\}$ be a Lévy process (meaning that $X$ has stationary and independent increments and is continuous in probability and that $X_0 = 0$), cadlag and centered, with moments of all orders. Note that Senguta [13] does not assume the stationarity of the increments in the definition of a Lévy process. Denote by $\sigma^2$ the variance of the Gaussian part of $X$ and by $\nu$ its Lévy measure. For background on all these notions, we refer to Sato [9]. The existence of moments of all orders of $X_t$ implies that the Lévy measure $\nu$ has moments of all orders $\geq 2$. Write

$$m_n = \int_{\mathbb{R}} x^n \, \nu(\mathrm{d}x) \qquad \text{for } n \geq 2.$$

Observe that, for $n \geq 3$, $m_n$ is the cumulant of $X_1$ of order $n$ and $m_2 + \sigma^2$ is the cumulant of order 2.

Following Nualart and Schoutens [8], we introduce the square-integrable martingales (and Lévy processes) called *Teugels martingales*, related to the powers of the jumps of the process:

$$Y_t^{(1)} = X_t,$$
$$Y_t^{(n)} = \sum_{0 < s \leq t} (\Delta X_t)^n - m_n t, \qquad n \geq 2,$$

where $\Delta X_s = X_s - X_{s-}$. The predictable quadratic variation of $Y^{(n)}$ is

$$\langle Y^{(n)} \rangle_t = C_n t,$$

where $C_n$ is a constant.

In the next theorem, we define a time–space harmonic polynomial related to the Lévy process $X$.

**Theorem 1.** *Let $X$ be a centered Lévy process with finite moments of all orders. Then, for every $n \geq 1$, the process*

$$M_t^{(n)} = \Gamma_n(X_t, -(m_2 + \sigma^2)t, -m_3 t, \ldots, -m_n t)$$

*is a martingale.*

**Proof.** Let $n \geq 2$. By Itô's formula, using (6) and the fact that $[X,X]_t^c = \sigma^2 t$, we have

$$M_t^{(n)} = n \int_0^t M_{s-}^{(n-1)} \, dX_s - \binom{n}{2}(m_2 + \sigma^2) \int_0^t M_s^{(n-2)} \, ds$$
$$- \sum_{j=3}^n \binom{n}{j} m_j \int_0^t M_s^{(n-j)} \, ds + \frac{1}{2} n(n-1) \sigma^2 \int_0^t M_s^{(n-2)} \, ds$$
$$+ \sum_{0 < s \leq t} (\Gamma_n(X_{s-} + \Delta X_s, -(m_2 + \sigma^2)s, -m_3 s, \ldots, -m_n s) \qquad (9)$$
$$- \Gamma_n(X_{s-}, -(m_2 + \sigma^2)s, -m_3 s, \ldots, -m_n s)$$
$$- n \Delta X_s \, \Gamma_{n-1}(X_{s-}, -(m_2 + \sigma^2)s, -m_3 s, \ldots, -m_n s)).$$

Applying (7) and rearranging terms, we obtain

$$\Gamma_n(X_{s-} + \Delta X_s, -(m_2 + \sigma^2)s, -m_3 s, \ldots, -m_n s) = \sum_{j=0}^n \binom{n}{j} M_{s-}^{(n-j)} (\Delta X_s)^j.$$

Then, the jumps part given in the expression for $M_t^{(n)}$ is

$$\sum_{0<s\leq t}\sum_{j=2}^{n}\binom{n}{j}M_{s-}^{(n-j)}(\Delta X_s)^j = \sum_{j=2}^{n}\binom{n}{j}\int_0^t M_{s-}^{(n-j)}\,\mathrm{d}X_s^{(j)}$$

$$= \sum_{j=2}^{n}\binom{n}{j}\int_0^t M_{s-}^{(n-j)}(\mathrm{d}Y_s^{(j)}+m_j\,\mathrm{d}s).$$

Therefore, after some cancellations, (9) can be written as

$$M_t^{(n)} = \sum_{j=1}^{n}\binom{n}{j}\int_0^t M_{s-}^{(n-j)}\,\mathrm{d}Y_s^{(j)}. \tag{10}$$

Finally, since, for every $t\geq 0$, we have $\mathbb{E}[\int_0^t M_{s-}^{(k)\,2}\,\mathrm{d}\langle Y^{(j)}\rangle_s]<\infty$ for $k\geq 0$ and $j\geq 1$, all the stochastic integrals on the right-hand side of (10) are martingales. □

**Corollary 2.** *Let $X$ be a centered Lévy process with finite moments of all orders. Define the polynomials*

$$Q_n(x,t) = \Gamma_n(x,-(m_2+\sigma^2)t,-m_3 t,\ldots,-m_n t), \qquad n\geq 1,$$

*and $Q_0(x,t)=1$. Then,*

(a) $Q_n(x,t)$, $n\geq 1$, *are time–space harmonic polynomials relative to $X$;*
(b) $Q_n(x,t)$, $n\geq 1$, *satisfy the recurrence relation*

$$Q_{n+1}(x,t) = xQ_n(x,t) - n(m_2+\sigma^2)tQ_{n-1}(x,t)$$
$$- \sum_{j=2}^{n}\binom{n}{j}m_{j+1}tQ_{n-j}(x,t);$$

(c) *an explicit expression for $M_t^{(n)}$ as a polynomial in $X_t$ is*

$$M_t^{(n)} = Q_n(X_t,t)$$
$$= \sum_{j=0}^{n}\binom{n}{j}\Gamma_{n-j}(0,-(m_2+\sigma^2)t,-m_3 t,\ldots,-m_{n-j}t)\,X_t^j, \qquad n\geq 1; \tag{11}$$

(d) *if $X$ is non-degenerate, then the family $\{Q_n(x,t),\ n\geq 1\}$ spans the set of all space–time harmonic polynomials relative to $X$.*

**Proof.** The recurrence in (b) is deduced from (5). Property (c) is obtained from (8). Finally, for (d), note that, from the properties stated in Section 2, the family $\{Q_n(x,t),\ n\geq 1\}$ satisfies the properties (i) to (iv) of Sengupta [13], page 953. On the other hand, if $X$ is non-degenerate, then the support of $X_t$ is unbounded (see [9], Theorem 24.3). Therefore the proof of Theorem 6 from Sengupta [13] can be applied to this family. □

**Remark 3.** 1. Assume that $X$ is non-degenerate. Let $X_t^n + p_{n-1}(t)X_t^{n-1} + \cdots + p_0(t)$ be a time–space harmonic polynomial relative to $X$, where $p_j(t)$ are ordinary polynomials in $t$. Formula (2.2) of Goswami and Sengupta [3] can be transferred to our context. Then, the coefficients of the polynomial must satisfy

$$p_i(s) = \sum_{j=i}^{n} \binom{j}{i} p_j(t) \mu_{j-i}(t-s), \qquad 0 \leq i \leq n, \ 0 \leq s \leq t,$$

where $\mu_r(t) = \mathbb{E}[X_t^r]$. Applying this property to $Q_n(x,t)$, from (11) and the fact that $\mu_r(t) = \Gamma_r(0, (m_2 + \sigma^2)t, m_3 t, \ldots, m_r t)$, it follows that

$$\Gamma_k(0, -(m_2 + \sigma^2)s, -m_3 s, \ldots, -m_k s)$$
$$= \sum_{\ell=0}^{k} \binom{k}{\ell} \Gamma_{k-\ell}(0, -(m_2 + \sigma^2)t, -m_3 t, \ldots, -m_{k-\ell} t)$$
$$\times \Gamma_\ell(0, (m_2 + \sigma^2)(t-s), m_3(t-s), \ldots, m_\ell(t-s)), \qquad 0 \leq s \leq t.$$

2. When the Lévy process $X$ has only a finite number of finite moments, say up to order $k$, then Theorem 1 is still true up to the polynomial of degree $[k/2] + 1$, where $[a]$ denotes the integer part of $a$.

## 4. Examples

### 4.1. Brownian motion

Let $\{W_t, t \geq\}$ be a standard Brownian motion and consider $X_t = W_t$. Then, $\sigma = 1$ and $\nu = 0$, so $m_n = 0$ for all $n \geq 2$. The time–space harmonic polynomials are

$$Q_n(x,t) = \Gamma_n(x, -t, 0, \ldots, 0),$$

which satisfy the recurrence relation

$$Q_{n+1}(x,t) = xQ_n(x,t) - ntQ_{n-1}(x,t), \qquad n \geq 1.$$

It is then clear that they coincide with the monic Hermite polynomials $\widetilde{H}_n(x,t)$ (we follow the notation of Schoutens [10], Chapter 5) defined by

$$\exp\{ux - tu^2/2\} = \sum_{n=0}^{\infty} \widetilde{H}_n(x,t) \frac{u^n}{n!}$$

or, equivalently,

$$\widetilde{H}_n(x,t) = \left(\frac{t}{2}\right)^{n/2} H_n(x/\sqrt{2t}),$$

where $H_n(t)$ is the classical Hermite polynomial given by

$$\exp\{2ux - u^2\} = \sum_{n=0}^{\infty} H_n(x)\frac{u^n}{n!}.$$

### 4.2. Poisson process

Let $\{N_t, t \geq 0\}$ be a Poisson process of intensity 1 and $X_t = N_t - t$ the compensated process. The Lévy measure of $X$ is a Dirac delta measure at point 1, hence $m_n = 1$, $n \geq 2$. Then,

$$Q_n(x,t) = \Gamma_n(x, -t, \ldots, -t),$$

and the recurrence relation is

$$Q_{n+1}(x,t) = xQ_n(x,t) - t\sum_{j=1}^{n}\binom{n}{j}Q_{n-j}(x,t).$$

Now, consider the monic Charlier polynomials $\widetilde{C}_n(x,t)$ (also with Schoutens' notation [10], Chapter 5), defined by

$$e^{-ut}(1+u)^x = \sum_{n=0}^{\infty} \widetilde{C}_n(x,t)\frac{u^n}{n!}.$$

It is well known that $\widetilde{C}_n(N_t, t) = \widetilde{C}_n(X_t + t, t)$ is a martingale. Write

$$\overline{C}_n(x,t) = \frac{1}{n!}\widetilde{C}_n(x+t,t).$$

Then (see, e.g., [10], Theorem 8),

$$\int_0^t \overline{C}_n(X_{s-}, s)\,dX_s = \overline{C}_{n+1}(X_t, t). \tag{12}$$

Corollary 2(d) guarantees that the polynomial $\overline{C}_n(x,t)$ is a linear combination of the polynomials $Q_j(x,t)$. To specify the relationship between both families of polynomials, we shall recursively construct the numbers $\lambda_1^{(n)}, \ldots, \lambda_n^{(n)}$ such that

$$Q_n(x,t) = \sum_{j=1}^{n} \lambda_j^{(n)} \overline{C}_j(x,t).$$

First,
$$Q_1(x,t) = x = \overline{C}_1(x,t),$$
so $\lambda_1^{(1)} = 1$. Assume that $\lambda_1^{(1)}, \lambda_1^{(2)}, \lambda_2^{(2)}, \ldots, \lambda_1^{(n-1)}, \ldots, \lambda_{n-1}^{(n-1)}$ are known. Then, from (10),
$$Q_n(X_t,t) = \sum_{j=1}^{n} \binom{n}{j} \int_0^t Q_{n-j}(X_{s-},s) \, \mathrm{d}Y_s^{(j)}.$$

For the compensated Poisson process, all the Teugels martingales are equal, that is $Y_t^{(n)} = X_t$, and using (12), we have
$$Q_n(X_t,t) = \sum_{j=1}^{n-1} \binom{n}{j} \sum_{k=1}^{n-j} \lambda_k^{(n-j)} \overline{C}_{k+1}(X_t,t) + \overline{C}_1(X_t,t)$$
$$= \sum_{k=1}^{n-1} \overline{C}_{k+1}(X_t,t) \sum_{j=1}^{n-k} \binom{n}{j} \lambda_k^{(n-j)} + \overline{C}_1(X_t,t).$$

It follows that $\lambda_1^{(n)} = 1$ and
$$\lambda_{k+1}^{(n)} = \sum_{j=k}^{n-1} \binom{n}{j} \lambda_k^{(j)}, \qquad k=1,\ldots,n-1.$$

### 4.3. Sum of two independent Lévy processes

Let $\{Y_t, t \geq 0\}$ and $\{Z_t, t \geq 0\}$ be two independent centered Lévy processes with moments of all orders, and put $X = Y + Z$. Since the cumulant of order $j$ of a sum of two independent random variables is the sum of the respectives cumulants of that order, the following convolution-type formula is deduced from (4):
$$\Gamma_n(y_1+z_1,\ldots,y_n+z_n) = \sum_{j=0}^{n} \binom{n}{j} \Gamma_j(y_1,\ldots,y_j) \Gamma_{n-j}(z_1,\ldots,z_{n-j}).$$

Hence,
$$Q_n^X(x,t) = \sum_{j=0}^{n} \binom{n}{j} Q_j^Y(y,t) Q_{n-j}^Z(z,t),$$
where $x = y + z$.

This result covers, for example, the process $X_t = W_t + N_t - t$, where $\{W_t, t \geq 0\}$ is a standard Brownian motion and $\{N_t, t \geq 0\}$ is a Poisson process of intensity 1, both processes being independent. The polynomials corresponding to $X$ are then convolutions between the Hermite polynomials and the polynomials studied in the previous example.

### 4.4. Gamma process

Let $\{G_t, t \geq 0\}$ be a Gamma process, that is, a Lévy process such that $G_t$ has distribution Gamma with mean $t$ and scale parameter equal to 1. Consider $X_t = G_t - t$. The Lévy measure is

$$\nu(\mathrm{d}x) = \frac{\mathrm{e}^{-x}}{x} \mathbf{1}_{(0,\infty)}(x)\,\mathrm{d}x$$

and the cumulants of $X_1$ of order $n \geq 2$ are

$$m_n = \int_0^\infty x^{n-1} \mathrm{e}^{-x}\,\mathrm{d}x = (n-1)!.$$

The polynomials are then

$$Q_n(x,t) = \Gamma_n(x, -t, \ldots, -(n-1)!\,t)$$

and the recurrence relation is

$$Q_{n+1}(x,t) = xQ_n(x,t) - n!\,t \sum_{j=0}^{n-1} \frac{1}{j!}\,Q_j(x,t).$$

We are going to prove, in agreement with the results of Schoutens and Teugels [11] (see also [10]), that $Q_n(x,t)$ is related to the Laguerre polynomials. First, we compute the generating function of the polynomials $Q_n(x,t)$. Applying the general formula (2) to this specific case, for $u$ in a neighborhood of zero of radius $\rho \in (0,1)$,

$$\sum_{n=0}^\infty Q_n(x,t)\frac{u^n}{n!} = \sum_{n=0}^\infty \Gamma_n(x, -t, \ldots, -(n-1)!\,t)\frac{u^n}{n!}$$

$$= \exp\left\{xu - t\sum_{n=2}^\infty \frac{u^n}{n}\right\}$$

$$= \exp\{u(x+t) + t\log(1-u)\}$$

$$= (1-u)^t \exp\{u(x+t)\}.$$

Second, using the relationship between Laguerre and Charlier polynomials, the generating function of the Laguerre polynomials can be written as (see [10], page 47)

$$(1-u)^t \exp\{ux\} = \sum_{n=0}^\infty (-1)^n L_n^{(t-n)}(x) u^n.$$

Then,

$$Q_n(x,t) = (-1)^n n!\,L_n^{(t-n)}(x+t).$$

Note that Schoutens [10] proves that $L_n^{(t-n)}(G_t)$ is a martingale, which is equivalent to the fact that $Q_n(X_t, t)$ is a martingale.

### 4.5. Compound Poisson process with lognormal jumps

The final example deals with a Lévy process with finite moments of all orders, but without moment generating function. Let $\{N_t, t \geq 0\}$ be a Poisson process of intensity 1 and $\{Y_n, n \geq 1\}$ a sequence of i.i.d. random variables with lognormal distribution of parameters $\mu = 0$ and $\sigma = 1$, independent of the Poisson process. Consider the centered compound Poisson process

$$X_t = \sum_{j=1}^{N_t} Y_j - e^{1/2} t,$$

with the convention that the sum is 0 when $N_t = 0$. As for all $\lambda > 0$, $\mathbb{E}[\exp\{\lambda |X_1|\}] = \infty$, it follows that $X_1$ does not have finite moment generating function. The Lévy measure is the law of $Y_1$ and the cumulants of $X_1$ of order $\geq 2$ are the moments of a lognormal law,

$$m_n = \mathbb{E}[Y_1^n] = e^{n^2/2}, \qquad n \geq 2.$$

The polynomials are then

$$Q_n(x, t) = \Gamma_n(x, -te^2, \ldots, -te^{n^2/2})$$

and the recurrence relation is

$$Q_{n+1}(x, t) = xQ_n(x, t) - t\sum_{j=1}^{n} \binom{n}{j} e^{(j+1)^2/2} Q_{n-j}(x, t).$$

## 5. Polynomials for a Lévy process with moment generating function

In this section, we study a general procedure (proposed by Sengupta [13]) for finding time–space harmonic polynomials related to a Lévy process with finite moment generating function in an open neighborhood of zero. We will prove that the resulting polynomials are the same as those obtained in Section 3.

Assume that $X_1$ has moment generating function in an open interval of the origin. Then,

$$\mathbb{E}[e^{uX_t}] = \exp\left\{\tfrac{1}{2}t\sigma^2 u^2 + t\int_{\mathbb{R}}(e^{ux} - 1 - ux)\,d\nu(x)\right\}$$

is analytic for $u \in (-\delta, \delta)$, for some $\delta > 0$, and it is known that (fixed $u$)

$$L_t(u) := \frac{\exp\{uX_t\}}{\mathbb{E}[\exp\{uX_t\}]} \\ = \exp\left\{uX_t - \frac{1}{2}t\sigma^2 u^2 - t\int_{\mathbb{R}}(e^{ux} - 1 - ux)\,d\nu(x)\right\} \quad (13)$$

is a martingale. On the other hand, for fixed $\omega \in \Omega$ and $t > 0$, as a function of $u$, $L_t(u)$ is also analytic in $(-\delta, \delta)$ and can be expanded as

$$L_t(u) = \sum_{n=0}^{\infty} R_n(X_t, t)\frac{u^n}{n!}. \quad (14)$$

Each coefficient $R_n(X_t, t)$ is a polynomial (in $X_t$ and $t$) that is a martingale (see the proof of a similar result in [10], page 46). Therefore, searching for a family of time–space harmonic polynomials is reduced to computing $R_n(x, t)$. This concludes the general procedure. However, we can go one step further since the condition that $X_1$ has moment generating function in an open interval of the origin is equivalent (see Remark 4 below) to the condition that for some $\lambda > 0$,

$$\int_{\{|x|>1\}} \exp\{\lambda|x|\}\,\nu(dx) < \infty.$$

Hence, the integral $\int_{\mathbb{R}}(e^{ux} - 1 - ux)\,d\nu(x)$ is analytic in a neighborhood of 0 and can be expanded as the power series $\sum_{n=2}^{\infty} \frac{u^n}{n!}m_n$. Substituting the integral by the power series in (13), we obtain

$$L_t(u) = \exp\left\{uX_t - \frac{1}{2}tu^2(\sigma^2 + m_2) - t\sum_{n=3}^{\infty}\frac{u^n}{n!}m_n\right\}. \quad (15)$$

Then, from (14), (15) and the exponential formula (2), we can identify the convenient $\mu_n$ and $k_n$ and obtain

$$R_n(X_t, t) = \Gamma_n(X_t, -t(\sigma^2 + m_2), -tm_3, \ldots, -tm_n) = Q_n(X_t, t).$$

**Remark 4.** The hypothesis that the random variable $X_1$ has finite generating function in some open neighborhood of zero means that there exists $\lambda > 0$ such that

$$\mathbb{E}[\exp\{\lambda|X_1|\}] < \infty. \quad (16)$$

The function $\exp\{\lambda|x|\}$ is submultiplicative and therefore (16) is equivalent to

$$\int_{\{|x|>1\}} \exp\{\lambda|x|\}\,\nu(dx) < \infty$$

(see [9], Theorem 25.3 and Proposition 25.4). It is worth noting that this last condition is equivalent to the hypothesis assumed by Nualart and Schoutens [8].

## 6. Time–space harmonic polynomials and orthogonal polynomials

Another interesting family of polynomials relative to a Lévy process $X$ with moments of all orders is defined in the following way. Consider the iterated integrals

$$P_t^{(0)} = 1, \qquad P_t^{(1)} = X_t, \ldots, P_t^{(n)} = \int_0^t P_{s-}^{(n-1)} \, dX_s. \tag{17}$$

These processes are related to the *variations* of $X$,

$$X_t^{(1)} = X_t, \qquad X_t^{(2)} = [X,X]_t, \qquad X_t^{(n)} = \sum_{s \leq t} (\Delta X_s)^n, \qquad n \geq 3, \tag{18}$$

through the Kailath–Segall formula

$$P_t^{(n)} = \frac{1}{n}(P_t^{(n-1)} X_t^{(1)} - P_t^{(n-2)} X_t^{(2)} + \cdots + (-1)^{n+1} P_t^{(0)} X_t^{(n)}) \tag{19}$$

(see Segall and Kailath [12] or Meyer [6]).

It follows that $P_t^{(n)}$ is a polynomial in $X_t^{(1)}, \ldots, X_t^{(n)}$, is called the *Kailath–Segall polynomial* of order $n$. Denote this polynomial by $P_n(x_1, \ldots, x_n)$, that is,

$$P_t^{(n)} = P_n(X_t^{(1)}, \ldots, X_t^{(n)})$$

and set $P_0 = 1$. The translation of (19) to the polynomials $P_n(x_1, \ldots, x_n)$ produces the recurrence relation

$$P_n(x_1, \ldots, x_n) \\
= \frac{1}{n}(P_{n-1}(x_1, \ldots, x_{n-1})x_1 - P_{n-2}(x_1, \ldots, x_{n-2})x_2 + \cdots + (-1)^{n+1} P_0 x_n). \tag{20}$$

On the other hand, since the iterated integrals of different order of a Lévy process are orthogonal,

$$\mathbb{E}[P_t^{(n)} P_t^{(m)}] = \frac{1}{n!} C^n t^n \delta_{nm}, \tag{21}$$

orthogonality of the polynomials $P_n(X_1^{(1)}, \ldots, X_t^{(n)})$ follows. However, only when the underlying process is a Brownian motion or a Poisson process is $P_n(X_1^{(1)}, \ldots, X_t^{(n)})$ a polynomial of the form $X_t^n + a_1(t) X_t^{n-1} + \cdots + a_n(t)$, with $a_j(t)$ being ordinary polynomials in $t$ (see Solé and Utzet [14]). In light of (11), this is a major difference from time–space harmonic polynomials.

In general, both families of polynomials are related by the polynomial $\Gamma_n$. Specifically,

$$P_n(x_1,\ldots,x_n) = \frac{1}{n!}\,\Gamma_n(x_1,-x_2,2!\,x_3,\ldots,(-1)^{n-1}(n-1)!\,x_n)$$

because, by (5), this polynomial satisfies the recurrence relation (20).

## Acknowledgements

This research was supported by Grant BFM2003-00261 of the Ministerio de Educación y Ciencia, Spain.